\documentclass[11pt]{amsart}
\usepackage{graphicx, comment}
\usepackage{times}
\usepackage{amsmath}
\usepackage{amssymb}
\usepackage{latexsym}
\newtheorem{theorem}{Theorem}[section]
\newtheorem{lemma}[theorem]{Lemma}

\theoremstyle{definition}
\newtheorem{definition}[theorem]{Definition}

\theoremstyle{remark}

\numberwithin{equation}{section}

\newcommand{\A}{\mathcal{A}}

\newcommand{\ddt}{\frac{d}{dt}}

\def \l {2}

\begin{document}

\title[Dissipation anomaly for a dyadic model] {The vanishing viscosity
limit for a dyadic model}

\author{Alexey Cheskidov}
\address{Department of Mathematics, Statistics and Computer Science,
University of Illinois at Chicago,
851 S Morgan St, Chicago, IL 60607}
\email{\tt acheskid@math.uic.edu}

\author{Susan Friedlander}
\address{Department of Mathematics,
University of Southern California
3620 South Vermont Ave., KAP 108
Los Angeles, CA 90089}
\email{\tt susan@math.northwestern.edu}

\date{\today}

\begin{abstract}
A dyadic shell model for the Navier-Stokes equations is studied
in the context of turbulence. The model is an infinite nonlinearly
coupled system of ODEs.
It is proved that the unique fixed point is a global attractor, which
converges to the global attractor  of the inviscid system
as viscosity goes to zero.
This implies that
the average dissipation rate
for the viscous system converges to the anomalous dissipation rate
for the inviscid system (which is positive) as viscosity goes to zero.
This phenomenon is called the dissipation anomaly predicted by Kolmogorov's
theory for the actual Navier-Stokes equations.
\end{abstract}
\maketitle

\section{Introduction}
We consider the Navier-Stokes equations for the motion
of a three-dimentional incompressible viscous fluid:
\begin{equation} \label{NSE}
\begin{split}
\frac{\partial u}{\partial t} + (u \cdot \nabla)u &= - \nabla p -\nu \Delta u +f,\\
\nabla \cdot u &=0,
\end{split}
\end{equation}
Here $u$ denotes the velocity vector field, $p$ the pressure, $f$
an external force, and $\nu$ the viscosity coefficient. The role
of the nonlinear term in \eqref{NSE} is critically important in the theory
of turbulence where a basic principle is a cascade of energy
from large scales, through the so called inertial scales, to very
small dissipative scales. Transfer of energy through these
scales is achieved via nonlinear interactions between the modes
in the Fourier space. This subject is the topic of extensive study
in the experimental, numerical, and analytical literature (see, for
example, Frisch~\cite{Fr}, Eyink and Sreenivasan \cite{ES}). Important
seminal work in the modern theory of turbulence was performed
by Kolmogorov and his school in the mid 20th century. However,
rigorous mathematical proofs of Kolmogorov's laws remain to
be obtained.

Kolmogorov predicted that the energy cascade mechanism in fully
developed three-dimensional turbulence produces a striking
phenomenon, namely the persistence of non-vanishing energy dissipation
in the limit of vanishing viscosity. This behavior, called
the ``dissipation anomaly'', can be described as follows
\[
\lim_{T \to \infty} \frac{1}{T} \int_0^T \nu \|u^\nu(t)\|_{H^1} \, dt \to \epsilon_\mathrm{d}>0,
\]
as $\nu \to 0$, where $u^\nu(t)$ is a solution to \eqref{NSE}.

Onsager conjectured that sufficiently rough solutions to the Euler
equations (i.e., \eqref{NSE} with $\nu=0$) can exhibit turbulent,
or anomalous dissipation. More precisely, if the H\"{o}lder exponent $h$
of the velocity is greater than $1/3$, then energy is conserved, however,
this ceases to be true if $h\leq 1/3$. For recent results concerning
Onsager's conjecture see, for example, Eyink \cite{E}, Constantin, E, and Titi \cite{cet}, Duchon and Robert \cite{DR}, Cheskidov, Constantin, Friedlander,
 and Shvydkoy \cite{CCFS}.

Partially because of the difficulty of proving mathematically
rigorous results in turbulence theory, a number of ``toy'' models that
preserve some features of the nonlinearity of the fluid equations have
been proposed and studied by physicists and mathematicians. These
include the so called shell models of the energy cascade, where the
nonlinearity of the 3D NSE  is simplified by considering only local interactions
between scales.
In this paper we study one of the original shell models introduced
in the context of oceanography by Desnyanskiy
and Novikov \cite{DN}.
This model, referred to as the dyadic model,
can be written as the following
infinite system of coupled ordinary differential equations:
\begin{equation} \label{e:dyadic}
\ddt a_j + \nu 2^{2j}a_j - 2^{c(j-1)}a_{j-1}^2 + 2^{cj}
a_j a_{j+1}=f_j, \qquad j=0,1,2,\dots,
\end{equation}
where $a_{-1}=0$, $c$ is a positive parameter, and $\frac{1}{2}a_j^2$ represents the total energy in the frequencies of order $2^j$. For convenience we chose the force $f$
so that $f_0>0$ and $f_j=0$ for all $j>0$.

This model has been analytically studied by Katz and Pavlovic \cite{KP}, 
Cheskidov \cite{C4}. Onsager's conjecture for the inviscid model (i.e.,
\eqref{e:dyadic} with $\nu=0$) was proved in Cheskidov, Friedlander,
and Pavlovic \cite{CFP1, CFP2}, where it was shown that the inviscid system
exhibits anomalous dissipation and the unique fixed point is a global
attractor.


Consider \eqref{NSE} on the whole space $\Omega = \mathbb{R}^3$ and
define $S_j u$ as follows:
\[
S_j u = u * \mathcal{F}^{-1} (\psi(\cdot 2^{-j})),
\]
where $\psi(\xi)$ is a smooth nonnegative function supported in the
ball of radius one centered at the origin and such that $\psi(\xi) =1$
for $\xi \leq 1/2$, and $\mathcal{F}$ is the Fourier transform. Then we define the energy flux due to nonlinear
interactions through the sphere of radius $2^j$ (see \cite{CCFS}) as
\[
\Pi_j = -\int_{\mathbb{R}^3} u \cdot \nabla S_j^2 u \cdot u \, dx.
\]
Using the test function $S^2_q u$ in the weak formulation of the
Navier-Stokes equations we obtain
\begin{equation} \label{e:flux1}
\frac{1}{2}\frac{d}{dt} S^2_j u = -\Pi_j - \nu \|\nabla S^2_j u\|_2.
\end{equation}
Recently, Cheskidov, Constantin, Friedlander, and Shvydkoy
\cite{CCFS} obtained the following new bounds on the nonlinear term 
in \eqref{NSE}:
\begin{equation} \label{e:local}
|\Pi_j| = \left| \int_{\mathbb{R}^3} u \cdot \nabla S_j^2 u \cdot u \, dx \right|
\lesssim  \sum_{i=-1}^\infty 2^
{-\frac{2}{3}|j-i|} 2^i \| u_i\|_3^3,
\end{equation}
where $u_j$ is a Littlewood-Paley piece of $u$ defined as
\[
u_j = S_{j+1}u - S_ju.
\]
This estimate employing the Littlewood-Paley decomposition
provides detailed information concerning the cascade of energy through frequency space. More precisely, it shows that the energy flux $\Pi_j$
through the sphere of radius $\kappa$ is controlled primarily by scales of order $\kappa$.

Recall that Bernstein's inequality can be stated as
\begin{equation} \label{e:Bernstein}
\|u_j\|_q \lesssim 2^{(3/p - 3/q)j} \|u_j\|_p, \qquad 1 \leq p \leq q.
\end{equation}
Now if we define
$a_j=\|u_j\|_2$, then using Bernstein's inequality with $p=2$ and $q=3$, we obtain
\begin{equation} \label{e:ineq}
a_j^3  \lesssim \|u_j\|_3^3 \lesssim 2^{3j/2} a_j^3,
\end{equation}
Motivated by \eqref{e:local} and \eqref{e:ineq}, we model the flux in the following way:
\begin{equation} \label{e:flux}
\Pi_j = 2^{cj}a_j^2 a_{j+1},
\end{equation}
where the scaling parameter $c \in [1,5/2]$. Here the bounds for the scaling
parameter $c$ are determined by Bernstein's inequality \eqref{e:Bernstein},
with the upper bound corresponding to saturation of the inequality.
In a turbulent flow it is expected that the degree of such saturation
could vary giving a rise to a phenomenon known as intermittency.
Motivated by \eqref{e:flux1} we write
\begin{equation} \label{j}
\frac{1}{2}\frac{d}{dt} \left(\sum_{i=0}^j a_i\right)^2 = -\Pi_j  - \nu
\sum_{i=0}^j 2^i a_i^2
\end{equation}
and
\begin{equation} \label{j-1}
\frac{1}{2}\frac{d}{dt} \left(\sum_{i=0}^{j-1} a_i\right)^2 = -\Pi_{j-1}  - \nu
\sum_{i=0}^{j-1} 2^i a_i^2.
\end{equation}
Subtracting \eqref{j-1} from \eqref{j} gives
\begin{equation}
\frac{1}{2}\frac{d}{dt} a_j^2 = \Pi_{j-1}-\Pi_j  - \nu 2^j a_j^2,
\end{equation}
which from the definition of the flux results
in the model system \eqref{e:dyadic}.

In this context, the
energy $E$ and the Sobolev norms are defined as
\[ 
E=\frac{1}{2}|a|^2= \frac{1}{2} \sum_{j=0}^\infty a_j^2, \qquad
\|a\|^2_{H^s} = \sum_{j=0}^\infty 2^{2js} a_j^2.
\] 

In this present paper we prove the following results for the viscous model
\eqref{e:dyadic}. In Section 2 we study a steady state $\alpha$, which
has a property that $\alpha_j$ is monotonic. This property is proved
for $3/2<c\leq 5/2$. Hence the rest of the results in this paper are valid only
for this range. In Section 3 we study the long-time behavior of solutions
to \eqref{e:dyadic} with $c>3/2$ and initial data $a(0) \in l^2$,
$a_j(0) \geq 0$ for all $j \geq 0$, and prove that the fixed point $\alpha$
is a global attractor. Moreover, $\alpha$ converges to the fixed point
of the inviscid system (which is a global attractor of the inviscid system)
as $\nu\to 0$. This allows us to conclude that the average dissipation rate
for the viscous system converges to the anomalous dissipation rate
for the inviscid system (which is positive) as $\nu \to 0$.

\subsection*{Acknowledgments}
A.C. was partially supported by NSF grant numbers DMS 0807827.
S.F. was partially supported by NSF grant numbers
DMS 0803268 and DMS 0503768.

\section{The fixed point}
In this section we study steady state solutions $\alpha$
to \eqref{e:dyadic}. We rescale the variables by
\[
A_j=\l^{cj/3} \l^{-c/6} f_0^{-1/2} \alpha_j
\]
to obtain the system of equations for steady states
\begin{equation} \label{e:steady}
\begin{split}
A_{j-1}^2 -A_j A_{j+1} &= \mu \l^{\beta j} A_j, \qquad j=1,2, \dots \\
-A_0 A_{1} &= \mu A_0 -1,
\end{split}
\end{equation}
where $\mu=\nu\l^{-c/6}f_0^{-1/2}$ and $\beta=2(1-c/3)$. For the inviscid
model studied in \cite{CFP1, CFP2}, there is a unique fixed point with
an explicit expression, namely $\{A_j=1\}$. In the case $\nu>0$
we cannot solve \eqref{e:steady} explicitly. However, we study
fixed points using the monotonicity property (Theorem~\ref{t:monoton}),
which we prove  following the analysis given by Heywood \cite{H}. He treated
the particular case of \eqref{e:steady} where $\beta$ was set to $0$.
In Section 3 we will show that the fixed point is a global attractor, i.e.,
it is unique.

We are only interested in finite energy solutions $a\in l^2$, which
translates to $A \in H^{-5/6}$.
The proof of the existence of a solution to \eqref{e:steady} follows
from standard Navier-Stokes techniques in which fixed point arguments
are used for truncations of the system. Standard arguments also show
that all these solutions are in $H^s$ for all $s>-5/6$. Moreover, we have
the following

\begin{lemma} \label{l:A_j-decay}
Let $A\in H^{-5/6}$ be a solution to \eqref{e:steady}. Then $A_j>0$ for all $j$ and
\[
\frac{A_{[J]+k}}{(A_{[J]+k-1})^2} \leq 2^{-\beta k}, \qquad k \geq 0,
\]
where $J$ is such that $\mu 2^{\beta J} = 1$.
\end{lemma}
\begin{proof}
Since $A \in H^s$ for all $s>0$, we have that $A_j \to 0$ as $j\to \infty$.
Then \eqref{e:steady} implies that
\[
\mu \sum_{j=0}^{\infty}\l^{\beta j} A_j^2 = A_0,
\]
and
\[
\mu \sum_{j=J}^{\infty}\l^{\beta j} A_j^2 = A_{J-1}^2A_J, \qquad J>0.
\]
Hence $A_j>0$ for all $j$. 

Now note that \eqref{e:steady} gives
\begin{equation}
\frac{A_{[J]+k}}{(A_{[J]+k-1})^2} \leq \frac{1}{\mu} 2^{-\beta([J]+k)}
\leq 2^{-\beta k}.
\end{equation}
%
\end{proof}

Note that since $A_j \to 0$ as $j\to \infty$, Lemma~\ref{l:A_j-decay}
implies that $A_j$ decays super-exponentially. This result does not
depend on the monotonicity
of the sequence $\{A_j\}$ and hence it holds in the whole range
$c\in[1, 5/2]$.

\begin{theorem} \label{t:monoton}
Every solution $A \in H^{-5/6}$ of \eqref{e:steady} with $c >3/2$ is monotonic, i.e.,
$A_{j-1} > A_{j}$ for all $j> 0$.
\end{theorem}
\begin{proof}
Let $h_j=A_j-A_{j-1}$. Then \eqref{e:steady} gives
\begin{equation} \label{e:h}
\begin{split}
h_{j+1}&=-h_j-\mu\l^{\beta j} - h_j A_{j-1}/A_j, \qquad j>0,\\
h_1&=-h_0-\mu+1/A_0.
\end{split}
\end{equation}
We prove that $h_j<0$ for all $j$ by contradiction. Assume that
$h_J \geq 0$ for some $J$. Then $h_{J+1} < -\mu \l^{\beta J}<0$, i.e.,
$A_{J+1}/A_J >1$. Then \eqref{e:h} implies that
\[
h_{J+2}>2\mu \l^{\beta J} - \mu\l^{\beta(J+1)} = \mu \l^{\beta J}(2-\l^\beta).
\]
Since $c>3/2$, we have that $\beta<1$. We conclude that $h_{J+2}>0$. Iterating this
process we obtain
\[
h_{J+2k-1} < -\mu \l^{\beta(J+2k-2)},
\]
and
\[
h_{J+2k} >0,
\]
for all $k>0$. Then
\[
A_{J+2k-2} > A_{J+2k-1} + \mu \l^{\beta(J+2k-2)},
\]
which contradicts the fact that $H^{s}$ norm of $A$ are finite for all $s>0$.
\end{proof}

\begin{lemma} \label{l:conv}
Let $A$ be a fixed point. Then
\[
\lim_{\mu \to 0} A_j = 1,
\]
for every $j\geq 0$.
\end{lemma}
\begin{proof}
The equations \eqref{e:steady} read
\begin{equation} \label{e:steady-mod}
\begin{split}
A_1 &= \frac{1}{A_0} - \mu,\\
A_{j+1} &= \frac{A_{j-1}^2}{A_{j}}-\mu \l^{\beta j}, \qquad j=0,1,2, \dots 
\end{split}
\end{equation}
We will proceed by induction. Suppose that $A_{j-1} \to 1$ as
$\mu \to 0$ for some $j\geq -1$. We will show that $A_{j} \to 1$
as $\mu \to 0$. Assume the contrary. Then we can pass to a subsequence
$\mu_n \to 0$ as $n\to \infty$, such that either
\[
\limsup_{\mu_n \to 0} A_j <1, \qquad \text{or} \qquad
\liminf_{\mu_n \to 0} A_j >1.
\]
First, assume that
$\limsup_{\mu_n \to 0} A_j <1$.
Then \eqref{e:steady-mod} implies that $\liminf_{\mu_n \to 0} A_{j+1} >1$,
which contradicts the monotonicity of $A$ (Theorem~\ref{t:monoton}).
Now assume that
$
\liminf_{\mu_n \to 0} A_j >1.
$
Then \eqref{e:steady-mod} implies that
\[
\limsup_{\mu_n \to 0} A_{j+1} <1, \qquad \text{and} \qquad
\liminf_{\mu_n \to 0} A_{j+2} >1,
\]
which again contradicts the monotonicity of $A$.
\end{proof}

The following property of a fixed point $A$ will be used in Section 3.

\begin{lemma} \label{l:prop}
Let $A$ be a fixed point. Then there exists $\gamma \in (0,1)$, such
that
\[
g_j(\mu):=\frac{A_{j+1}}{A_j+A_{j+1}^{1/2}A_{j+2}^{1/2}} < (1-\gamma) \l^{-\beta/2},
\]
for all $\mu>0$ and $j\geq 0$.
\end{lemma}
\begin{proof}
The monotonicity of $A$ implies
\begin{equation} \label{e:fj-start}
g_j(\mu) <\frac{A_{j+1}}{A_j+A_{j+2}}.
\end{equation}
From \eqref{e:steady} and monotonicity we also obtain
\begin{equation}
A_{j+1}^2-\mu\l^{\beta(j+1)}A_{j+2}=A_{j+2}A_{j+3}<A_{j+2}^2.
\end{equation}
Hence
\begin{equation} \label{e:longexp}
A_{j+2}> -{\textstyle \frac{1}{2}}\mu\l^{\beta(j+2)}+
{\textstyle \frac{1}{2}}\sqrt{\mu^2\l^{2\beta(j+2)}+4A_{j+1}^2}.
\end{equation}
We define
\[
y=\frac{A_{j+1}}{A_j}, \qquad z=\frac{\mu\l^{\beta(j+2)}}{2A_j}.
\]
From \eqref{e:steady} and the positivity of each $A_j$ we conclude
\[
A_{j+1} < \frac{A_j^2}{\mu \l^{\beta(j+1)}}.
\]
Hence,
\begin{equation} \label{e:constraint}
yz < 2^{\beta-1}.
\end{equation}
Substituting \eqref{e:longexp} into \eqref{e:fj-start} gives
\[
g_j(\mu)<\frac{y}{1-z+\sqrt{y^2+z^2}} 
\]
subject to constraints \eqref{e:constraint}, $0\leq y \leq 1$, and
$0\leq z < \infty$.
Hence
\[
g_j(\mu)<\frac{y^2}{y-2^{\beta-1}+\sqrt{y^4+2^{2(\beta-1)}}} =:h(y,\beta),
\]
Since $\frac{\partial h}{\partial y} >0$, $h$ attains its maximum at $y=1$.
Thus
\[
g_j(\mu)<\frac{1}{1-2^{\beta-1}+\sqrt{1+2^{2(\beta-1)}}} < (1-\gamma)2^{-\beta/2},
\]
provided $\beta <1$, i.e., provided $c>3/2$.
%
\end{proof}

\section{The global attractor}
In this section we study the long-time behavior of solutions to the time
dependent system 
We study the viscous dyadic model
\begin{equation} \label{model}
\ddt a_j -  2^{c(j-1)} a_{j-1}^2 + 2^{cj} a_j a_{j+1} + \nu 2^{j}a_j =f_j, \qquad j=0,1,2\dots,
\end{equation}
where $a_{-1}=0$. Here $c \in (3/2, 5/2]$, and $\nu >0$
is the viscosity. The initial data is assumed to be $a(0) \in l^2$,
$a_j(0) \geq 0$ for all $j$. 

\begin{definition}
A solution of \eqref{model} on $[T,\infty)$ (or $(-\infty, \infty)$, if
$T=-\infty$) of \eqref{model} is an $l^2$-valued
function $a(t)$ defined for $t \in [T, \infty)$, such that $a_j \in C^1([T,\infty))$
and $a_j(t)$ satisfies \eqref{model} for all $j$.
\end{definition}

Note that if $a(t)$ is a solution on $[T,\infty)$, then automatically
$a_j \in C^{\infty}([T,\infty))$. 
The following theorems were proved in \cite{C}.


\begin{theorem} \label{t:pos}
For every $a^0 \in l^2$ with $a^0_j\geq 0$ there exists
a solution of \eqref{model} with $a(0)=a^0$. Moreover,
$a_j(t) \geq 0$ for all $t>0$.
\end{theorem}

\begin{theorem} \label{t:iquality}
Let $a(t)$ be a solution to  \eqref{model} with
$a_j(0)\geq 0$. Then $a(t)$ satisfies the energy inequality
\begin{equation} \label{ee}
|a(t)|^2 + 2\nu \int_{t_0}^t \|a(\tau)\|_{H^1}^2 \, d\tau \leq
|a(t_0)|^2 + 2\int_{t_0}^t (f, a(\tau)) \, d\tau,
\end{equation}
for all $0 \leq t_0 \leq t$.
\end{theorem}

We write a
solution $a(t)$ to \eqref{model} in the form
\begin{equation} \label{e:resc-a}
a_j(t)=\alpha_j + b_j(t),
\end{equation}
where $\alpha$ is a fixed point whose properties were exhibited
in Section~2. We now show that this fixed point is
the exponential global attractor. In particular, the fixed
point is unique.

\begin{theorem}
Let $\alpha \in l^2$ be a fixed point of \eqref{model} for $c \in (3/2,5/2]$
and $a(t)$ be
a solution with $a(0) \in l^2$ and $a_j(0) \geq 0$ for all $j$. Then
\begin{equation}
|b(t)|^2 \leq  |b(0)|^2 e^{-2\gamma \nu t}.
\end{equation}
\end{theorem}
\begin{proof}
As before we write
\[
\alpha_j = 2^{-cj/3}2^{c/6}f_0^{1/2}A_j.
\]
Now let $B_j(t):=\l^{cj/3}b_j(t)$. Then
\begin{equation} \label{e:resc-a}
a_j(t)=\alpha_j + \l^{-cj/3}B_j(t).
\end{equation}
Then the system \eqref{model} reduces to
\begin{equation} \label{e:B}
2^{-2cj/3}f_0^{-1/2}2^{c/6} \frac{dB_j}{dt} = 2A_{j-1}B_{j-1} +B_{j-1}^2 -A_jB_{j+1}
-A_{j+1}B_j -B_jB_{j+1} - \mu 2^{\beta j} B_j,
\end{equation}
where $j\geq 0$ and $A_{-1}=B_{-1}=0$.
Following the procedures in the inviscid case given in \cite{CFP2} we
multiply \eqref{e:B} by $B_j$ and sum to obtain
\begin{equation}
f_0^{-1/2}2^{c/6}\frac{1}{2} \frac{d}{dt} \sum_{j=0}^k 2^{-2cj/3} B_j^2 = \sum_{j=0}^k(2A_{j-1}B_{j-1}B_j
-A_jB_{j+1}B_j - A_{j+1}B_j^2 - \mu 2^{\beta j} B_j^2 ) - B_k^2B_{k+1}.
\end{equation}
Since $A_k+B_k(t) \geq 0$ for all $t\geq 0$ and $\lim_{k\to \infty}A_k =0$,
we have that $\liminf_{k\to \infty} B_{k+1}(t) \geq 0$ for all $t\geq 0$.
Then due to the fact that $A_j$ decreases super-exponentially (see 
Lemma~\ref{l:A_j-decay}), and $\sum_{j=0}^\infty 2^{\beta j} B_j^2$
is integrable, we can use the dominated convergence theorem to
obtain
\begin{multline}
f_0^{-1/2}2^{c/6}\sum_{j=0}^\infty \l^{-2cj/3}B_j(t)^2 - f_0^{-1/2}2^{c/6}\sum_{j=0}^\infty \l^{-2cj/3}B_j(0)^2\\
\leq 2\int_{0}^t \sum_{j=0}^\infty \left[
2A_{j-1}B_{j-1}B_j
-A_jB_{j+1}B_j - A_{j+1}B_j^2 - \mu 2^{\beta j} B_j^2
\right] \, d\tau,
\end{multline}
where $j \geq 0$ and $A_{-1}=B_{-1}=0$.
Then we have
\begin{equation}
\begin{split}
\sum_{j=0}^\infty&\left[
2A_{j-1}B_{j-1}B_j
-A_jB_{j+1}B_j - A_{j+1}B_j^2 - \nu 2^{\beta j} B_j^2 \right] \\
\leq &-\frac{1}{2} \sum_{j=0}^\infty(A_{j+1}^{1/2}B_j - A_{j+2}^{1/2}B_{j+1})^2
+\sum_{j=0}^\infty(A_{j} - A_{j+1}^{1/2}A_{j+2}^{1/2})B_j B_{j+1}\\
&- \mu \sum_{j=0}^\infty \l^{\beta j}B_j^2\\
\leq &-\frac{1}{2} \sum_{j=0}^\infty(A_{j+1}^{1/2}B_j - A_{j+2}^{1/2}B_{j+1})^2
+
 \sum_{j=0}^\infty \frac{\mu A_{j+1}}{A_j+A_{j+1}^{1/2}A_{j+2}^{1/2}}
\l^{\beta(j+1)}B_jB_{j+1}\\
&- \mu \sum_{j=0}^\infty \l^{\beta j}B_j^2,
\end{split}
\end{equation}
where we used equation~\eqref{e:steady} in the last inequality.
Now using Cauchy-Schwarz inequality and Lemma~\ref{l:prop} we obtain
\begin{equation}
\begin{split}
\sum_{j=0}^\infty&\left[
2A_{j-1}B_{j-1}B_j
-A_jB_{j+1}B_j - A_{j+1}B_j^2 - \mu 2^{\beta j} B_j^2 \right] 
\leq  -\gamma\mu \sum_{j=0}^\infty \l^{\beta j}B_j^2.
\end{split}
\end{equation}
Therefore
\begin{equation}
\sum_{j=0}^\infty \l^{-2cj/3}B_j(t)^2 - \sum_{j=0}^\infty \l^{-2jc/3}B_j(0)^2
\leq  2f_0^{1/2}2^{-c/6} \gamma\mu \int_0^t \sum_{j=0}^\infty \l^{\beta j}B_j(\tau)^2 \, d\tau
\end{equation}
Hence for $b_j=\l^{-cj/3}B_j$ we have
\begin{equation}
\begin{split}
\sum_{j=0}^\infty b_j(t)^2 - \sum_{j=0}^\infty b_j(0)^2
&\leq -2\gamma\nu \int_0^t \sum_{j=0}^\infty 2^{2j} b_j(\tau)^2 \, d\tau\\
&\leq -2\gamma\nu \int_0^t \sum_{j=0}^\infty b_j(\tau)^2 \, d\tau,
\end{split}
\end{equation}
which implies that
\begin{equation}
|b(t)|^2 \leq  |b(0)|^2 e^{-2\gamma\nu t}.
\end{equation}
\end{proof}

\section{Dissipation Anomaly}

Here we study the energy dissipation in the limit of
vanishing viscosity. For convenience, solutions to the
dyadic model with viscosity $\nu \geq 0$ will be denoted by
$a^\nu(t)$ in this section. The fixed point (which is unique
in both viscous and inviscid cases) will be denoted by $\alpha^\nu$.
Now given a solution $a^0(t)$ to the inviscid dyadic model, we define its anomalous energy dissipation rate
as follows
\[
\epsilon_{a^0}(t) := (a^0(t),f) - \frac{1}{2} \frac{d}{dt} |a^0(t)|^2
\]
in the sense of distributions. Due to the energy inequality,
$\epsilon_{a^0} \geq 0$
and hence $\epsilon_{a^0}$ is a Borel measure.
The following theorem was proved in \cite{CFP}:
\begin{theorem} \label{t:dissipation}
Let $a^0(t)$ be a solution to the inviscid dyadic model on $[0,\infty)$. Then
\[
\lim_{T \to \infty} \frac{1}{T} \int_0^T d \epsilon_{a^0}(t) =: \epsilon_\mathrm{d}>0.
\]
\end{theorem}
Since the global attractor of the viscous model $\A^\nu$ converges to the
global attractor of the inviscid model $\A^0$ as $\nu \to 0$ (see Lemma~\ref{l:conv}), we have the following
\begin{theorem}
Let $a^\nu(t)$ be a solution to the viscous dyadic model on $[0,\infty)$.
Then
\[
\lim_{T \to \infty} \frac{1}{T} \int_0^T \nu \|a^\nu(t)\|_{H^1}^2 \, dt \to \epsilon_\mathrm{d}>0,
\]
as $\nu \to 0$.
\end{theorem}
\begin{proof}
Due to the energy inequality we have
\[
\frac{1}{2T}|a^\nu(t)|^2 - \frac{1}{2T}|a^\nu(t+T)|^2 \leq -\nu \frac{1}{T}
\int_t^{t+T} \|a^\nu(s)\|_{H^1}^2 \, ds + \frac{1}{T}\int_t^{t+T} (a^\nu(s),f) \, ds.
\]
Hence,
\[
\limsup_{T \to \infty} \frac{1}{T} \int_0^T \nu \|a^\nu(t)\|_{H^1}^2 \, dt \leq
\lim_{T \to \infty} \frac{1}{T}\int_0^{T} (a^\nu(s),f) \, ds =
(\alpha^\nu, f).
\]
On the other hand, note that the fixed point $\alpha^\nu$
(since it is a regular solution) satisfies the energy equality
\[
\nu \|\alpha^\nu\|_{H^1}^2 = (\alpha^\nu, f).
\]
Now for any $\delta>0$, there exists $N$, such that
\[
\nu \sum_{j=0}^N 2^{2j} (\alpha_j^\nu)^2 \geq \nu \|\alpha^\nu\|_{H^1}^2 - \delta.
\]
Since $a^\nu(t) \to \alpha^\nu$ in $l^2$,  we have
\[
\liminf_{T \to \infty} \frac{1}{T} \int_0^T \nu  \sum_{j=0}^N 2^{2j} a^\nu(t)_j^2\, dt \geq \nu \sum_{j=0}^N 2^{2j} (\alpha_j^\nu)^2 \geq \nu \|\alpha^\nu\|_{H^1}^2 - \delta.
\]
Therefore,
\[
\begin{split}
\liminf_{T \to \infty} \frac{1}{T} \int_0^T \nu \|a^\nu(t)\|_{H^1}^2 \, dt \geq
\nu\|\alpha^\nu\|_{H^1}^2 = (\alpha^\nu, f).
\end{split}
\]

Since $(\alpha^\nu, f) \to (\alpha^0, f)$ as $\nu \to 0$, we obtain
\[
\lim_{T \to \infty} \frac{1}{T} \int_0^T \nu \|a^\nu(t)\|_{H^1}^2 \, dt = (\alpha^0, f) = \alpha^0_0 f_0=: \epsilon_{\mathrm{d}} >0.
\]
\end{proof}

\section{Dissipation length scale}

Kolmogorov's theory of turbulence predicts that the energy density 
in the inertial range is
\begin{equation}\label{en-K}
\mathcal{E}(\kappa) \sim \epsilon_{\mathrm{d}}^{2/3} \kappa^{-5/3},
\end{equation}
followed by a rapid decay after the dissipation wave number
\begin{equation} \label{dis-K}
\kappa_d \sim \left(\frac{\epsilon_{\mathrm{d}}}{\nu^3}\right)^{\frac{1}{4}}.
\end{equation}
For the dyadic model the energy density is defined as $\mathcal{E}(2^j) = a_j^2 2^{-j}$.
In the inviscid case one can easily check that the energy density
for the fixed point (which is a global attractor) is
\begin{equation} \label{en-d}
\mathcal{E} (\kappa) \sim \epsilon_{\mathrm{d}}^{2/3} \kappa^{-2c/3-1}.
\end{equation}
By Lemma~\ref{l:conv} the energy density on the global attractor
of the dyadic system with small positive viscosity is close to \eqref{en-d}
in the inertial range. Moreover, Lemma~\ref{l:A_j-decay} can
be used to determine the dissipation wavenumber for the model
\begin{equation} \label{dis-d}
\kappa_d := 2^J = \left(\frac{f_0^{3/2}}{\nu^3}\right)^{\frac{1}{4}\cdot\frac{2}{3-c}}
\sim \left( \frac{\epsilon_{\mathrm{d}}}{\nu^3}\right)^{\frac{1}{4}\cdot\frac{2}{3-c}}.
\end{equation}
Lemma 2.1 is valid for all c in our range of interest, 
including $c = 1$. The inequality (2.2) implies that ultimately
the energy in each shell decays very rapidly with increasing $j$.
When we invoke monotonicity of the sequence $A_j$, (2.2) ensures that
this very rapid decay occurs at the dissipation wave number scale
given by (5.4). However for technical reasons we require $c > 3/2$
in order to prove monotonicity.

As we discussed in Section 1, the appropriate range for $c$, that
arises from Bernstein's inequality applied to estimate the nonlinear term
in the Navier-Stokes equations, is $c\in[1,5/2]$. 
We note in the end point case $c=1$ the expressions
for the energy density and the dissipation wavenumber \eqref{en-d} and
\eqref{dis-d} coincide
with \eqref{en-K} and \eqref{dis-K}, i.e.,
\[
\mathcal{E} (\kappa) \sim \epsilon_{\mathrm{d}}^{2/3} \kappa^{-5/3},
\qquad \kappa_d \sim \left(\frac{\epsilon_{\mathrm{d}}}{\nu^3}\right)^{\frac{1}{4}}, \qquad \text{for} \quad c=1.
\]
At the end point value $c=5/2$, which corresponds to complete
saturation of Bernstein's (and Sobolev) inequalities, we have
\[
\mathcal{E} (\kappa) \sim \epsilon_{\mathrm{d}}^{2/3} \kappa^{-8/3},
\qquad \kappa_d \sim \frac{\epsilon_{\mathrm{d}}}{\nu^3}, \qquad
\text{for} \quad c=5/2.
\]



\begin{thebibliography}{99}

\bibitem{C4} A. Cheskidov, Blow-up in finite time for dyadic models of the Navier-Stokes equations,
{\em Trans. Amer. Math. Soc.} {\bf 360} (2008), 5101-5120.  

\bibitem{C} A. Cheskidov,
Global attractors of evolutionary systems, {\it Journal of dynamics and Differential Equations,} to appear, arXiv:math.DS/0609357.

\bibitem{CFP1} A. Cheskidov, S. Friedlander and N. Pavlovi\'{c}, An inviscid dyadic model of turbulence:
the fixed point and Onsager's conjecture,
{\it J. Math. Phys.} {\bf 48}, 065503 (2007).

\bibitem{CCFS}
A. Cheskidov, P. Constantin, S. Friedlander, and R. Shvydkoy,
Energy conservation and Onsager's conjecture for the Euler equations,
{\it Nonlinearity}, to appear
arXiv:0704.0759v1 [math.AP].

\bibitem{CFP2} A. Cheskidov, S. Friedlander and N. Pavlovi\'{c}, An inviscid  dyadic model of turbulence: the global attractor, {\it DCDS-A}, to appear.

\bibitem{cet} P. Constantin, W. E, E. Titi, Onsager's conjecture on the energy conservation for
solutions of Euler's equation, {\em Commun. Math. Phys.} {\bf 165} (1994), 207--209.

\bibitem{DN} V. N. Desnyansky and E. A. Novikov, The evolution of turbulence spectra to the
similarity regime, {\it lzv. Akad Nauk SSSR Fiz. Atmos. Okeana} {\bf{10}} (1974), 127--136. 


\bibitem{DR} J. Duchon and R. Robert, Inertial energy dissipation for weak solutions of
incompressible Euler and Navier-Stokes equations, {\em Nonlinearity} {\bf 13} (2000), 249--255.

\bibitem{E} G. L. Eyink, Energy dissipation without viscosity in ideal hydrodynamics. I.
Fourier analysis and local energy transfer, {\em Phys. D} {\bf 78} (1994), 222--240.

\bibitem{ES} G. L. Eyink and K. R. Sreenivasan, Onsager and the theory of hydrodynamic
turbulence, {\it Rev. Mod. Phys.} {\bf{78}} (2006).

\bibitem{Fr} U. Frisch, 
{\em Turbulence: The legacy of A.N. Kolmogorov}, 
Cambridge University Press, Cambridge (1995), The legacy of A.N. Kolmogorov.

\bibitem{H} J. G. Heywood, A curious phenomenon in a model problem, suggestive of the hydrodynamic inertial range and smallest scale of motion,
{\it J. math. fluid mech.} {\bf 5} (2003), 403--423.

\bibitem{KP} N. H. Katz and N. Pavlovi\'c, Finite time blow-up for a dyadic model
of the Euler equations, {\it AMS Tran.} {\bf 357} (2005), 695--708. 


\end{thebibliography}
\end{document}